\documentclass[oneside]{amsart}

\usepackage[letterpaper,body={12.6cm,20.5cm}, mag=1000]{geometry}
\usepackage{amssymb}
\usepackage{amsthm}
\usepackage{amscd}

\numberwithin{equation}{section}

\theoremstyle{plain}

\newtheorem{lemma}[equation]{Lemma}

\newtheorem{prop}[equation]{Proposition}
\newtheorem{step}{Step}
\newtheorem{thm}[equation]{Theorem}

\newtheorem*{thma}{Theorem A}
\newtheorem*{thmb}{Theorem B}

\theoremstyle{definition}

\newcommand{\dlabel}[1]{\ifmmode \text{\ttfamily \upshape [#1] } \else
{\ttfamily \upshape [#1] }\fi \label{#1}}

\newcommand{\C}{\operatorname{C} }

\newcommand{\Z}{\operatorname{Z} }

\newcommand{\Cl}{\operatorname{Cl} }
\newcommand{\cl}{\operatorname{cl} }

\newcommand{\I}{\operatorname{I} }

\newcommand{\gen}[1]{\left < #1 \right >}

\newcommand{\Aut}{\operatorname{Aut} }

\newcommand{\Hom}{\operatorname{Hom} }
\newcommand{\Inn}{\operatorname{Inn} }

\newcommand{\Autcent}{\operatorname{Autcent} }

\newcommand{\annd}{\quad \text{ and } \quad}

\begin{document}

\title{Class preserving automorphisms of finite $p$-groups}

\author{Manoj K.~Yadav}

\address{School of Mathematics, Harish-Chandra Research Institute \\
Chhatnag Road, Jhunsi, Allahabad - 211 019, INDIA}

\email{myadav@mri.ernet.in}
\thanks{Research supported by DST (SERC Division), the Govt. of INDIA}

\date{\today}

\begin{abstract}
We classify all finite $p$-groups $G$ for which $|\Aut_{c}(G)|$  attains its 
maximum value, where $\Aut_{c}(G)$ denotes the group of all class preserving 
automorphisms of $G$.
\end{abstract}

\maketitle

\section{Introduction}

Let $G$ be a finite $p$-group and $|G| = p^n$, where $p$ is a prime and $n$ is
a non-negative integer. For $x \in G$, $x^G$ denotes the conjugacy class of 
$x$ in $G$. By $\Aut(G)$ we denote the group of all automorphisms of $G$.
An automorphism $\alpha$ of $G$ is called \emph{class preserving} if 
$\alpha(x) \in x^G$ for all $x \in G$. The set of all class preserving
automorphisms of $G$,  denoted by $\Aut_{c}(G)$, is a normal subgroup of 
$\Aut(G)$.

In 1911, W. Burnside \cite{wB55} posed the following question: Does there exist
any  finite group $G$ such that $G$ has a non-inner class preserving 
automorphism?
In 1913,  Burnside \cite{wB13} himself gave an affirmative answer to this 
question.
He constructed a group $G$ of order $p^6$ isomorphic to the group $W$ 
consisting  of all $3 \times 3$ matrices
\[  M = \begin{pmatrix} 
          1 & 0 & 0 \\
          x & 1 & 0 \\
          z & y & 1
        \end{pmatrix}  \]
with $x,y,z$ in the field $\mathbb F_{p^2}$ of $p^2$ elements, where $p$ is an
odd prime.
For this group $G$,  $\Inn(G) < \Aut_{c}(G)$, where $\Inn(G)$ denotes the
group of all inner automorphisms of $G$.
He also proved that $\Aut_{c}(G)$ is an elementary abelian $p$-group of order
$p^8$. Throughout the paper this group $G$ is represented by the group $W$.

For the statement of our main theorem we need the following theorem, 
which is also of independent interest as it provides a very neat bound 
for $|\Aut_{c}(G)|$:
\begin{thma} 
Let $G$ be a non-trivial $p$-group having order $p^{n}$. Then
\begin{equation}
\label{gbineq} |\Aut_{c}(G)| \leq
   \begin{cases}
     p^{\frac{(n^{2}-4)}{4}},  &\text{if $n$ is even;}\\
     p^{\frac{(n^{2}-1)}{4}}, &\text{if $n$ is odd.}
   \end{cases}
\end{equation}
\end{thma}
We prove Theorem A in Section $5$ as Theorem \ref{thm4}.

One can easily notice that equality holds in \eqref{gbineq} for the group $W$. 
Motivated from this we have the following natural problem: 

\noindent{\bf Problem.} Classify all finite $p$-group $G$  
such that equality holds in \eqref{gbineq}.

 In this paper, which may be viewed as a continuation of \cite{wB13}, 
we solve this problem.
For the statement of our main theorem we also need the 
definition of  Camina $p$-groups. 
These are the finite $p$-groups $G$ such that each 
non-trivial coset $x\gamma_2(G)$ of the derived group $\gamma_2(G)$ is a single
conjugacy class $x^G$ in $G$. For the detailed definition of Camina groups one
can see Section $3$ below. Now we state our main theorem, which we prove 
in Section $5$ as Theorem \ref{thm6}.
 
\begin{thmb}
Let $G$ be a non-abelian finite $p$-group of order $p^n$. Then
equality holds in \eqref{gbineq} if and only if one of the following holds: 
\begin{subequations}
\begin{align}
& \text{$G$ is an extra-special $p$-group of order $p^3$;}\\
& \text{$G$ is  a group of nilpotency class $3$ and order $p^4$;}\\
& \text{$G$ is a Camina special $p$-group isoclinic to the group $W$ and
    $|G|=p^6$;}\\
& \text{$G$ is isoclinic to $R$ and $|G|=p^6$,}
\end{align}
\end{subequations}
where $R$ is defined in \eqref{r}.
\end{thmb}

\noindent{\bf Acknowledgements.} I thank Prof. Everett C. Dade 
for his valuable comments, suggestions and corrections and  
Prof. I. B. S. Passi for advising me to work on `bound for $|\Aut_c(G)|$'.

\section{Notation}

Our notation for objects associated with a finite multiplicative group $G$ is 
mostly standard.  We use $1$ to denote both the identity element of 
$G$ and the trivial subgroup $\{ 1 \}$ of $G$. By $\Aut(G)$, $\Aut_{c}(G)$ and
$\Inn(G)$, we denote the group of all automorphisms, the group of conjugacy 
class 
preserving automorphisms and the group of inner automorphisms of $G$
respectively. The abelian group of all homomorphisms from an abelian group 
$H$ to an abelian group $K$ is denoted by $\Hom(H,K)$.

We write $\gen x$ for the cyclic subgroup of $G$ generated by a given element 
$x \in G$. To say that some $H$ is a subset or a subgroup of $G$ we write 
$H \subseteq G$ or $H \leq G$ respectively. To indicate, in addition, that 
$H$ is properly contained in $G$, we write $H \subset G$, $H < G$ 
 respectively.  If $x,y \in G$, then $x^y$ denotes the conjugate element 
$y^{-1}xy \in G$ and $[x,y] = [x,y]_G$ denotes the commutator 
$x^{-1}y^{-1}xy = x^{-1}x^y \in G$. If $x \in G$, then $x^G$ 
denotes the $G$-conjugacy class of all $x^w$, for $w \in G$, and $[x,G]$ 
denotes the set of all $[x,w]$, for $w \in G$. Since $x^w = x[x,w]$, for all 
$w \in G$, we have $x^G = x[x,G]$. For $x \in G$, $\C_{H}(x)$ denotes the 
centralizer of $x$ in $H$, where $H \le G$. The center of $G$ will be 
denoted by $\Z(G)$.

By long-standing convention the expression $[K,H]$, for subgroups $K,H \le G$, 
denotes, not the set of all commutators $[x,y]$, for $x \in K$ and $y \in H$, 
but rather the subgroup $\gen{\, [x,y] \mid x \in K, y \in H\,}$ of $G$ 
generated by those commutators. This does not conflict with the previous 
notation when $K$ is the element or subgroup $1$ of $G$, since both the subset 
$[1,H]$ and the subgroup $[1,H]$ are equal to $\{ 1 \}$.  

We write the subgroups in the lower central series of $G$ as $\gamma_n(G)$, 
where $n$ runs over all strictly positive integers. They are defined 
inductively by
\begin{subequations} \label{LCSOfG}
\begin{align} 
\gamma_1(G) &= G \annd \label{LCSOne} \\
\gamma_{n+1}(G) &= [\gamma_n(G), G] \label{LCSNPlusOne}
\end{align}
\end{subequations}
for any integer $n \ge 1$. Note that $\gamma_2(G)$ is the derived group 
$[G,G]$ of $G$. 
Let $x_1, x_2,\cdots,x_k$ be $k$ elements of $G$, where $k \ge 2$. 
The commutator of  $x_1$ and $x_2$ has been defined to be 
$[x_1,x_2] = x_1^{-1}x_2^{-1}x_1x_2$.
Now we define a higher commutator of $x_1, x_2,\cdots,x_k$ inductively as 
\[[x_1,x_2,\cdots,x_k] = [[x_{1},\cdots,x_{k-1}],x_k].\]
We have already defined a commutator subgroup $[H_1,H_2]$ of any
two subgroups $H_1$ and $H_2$ of $G$.
A higher commutator subgroup of any $k$ subgroups $H_1, \cdots, H_k$ of $G$ is
defined inductively by the formula
\[[H_1,H_2,\cdots,H_k] = [[H_1,\cdots,H_{k-1}],H_k],\]
where $k \ge 2$.
We will be using the commutator identities
\[[x,yz] = [x,z][x,y]^z = [x,z][x,y][x,y,z]\]
and
\[[xy,z] = [x,z]^y[y,z] = [x,z][x,z,y][y,z]\]
many times without any reference.

\section{Camina groups of class $3$}

Let $G$ be a finite group and $1 \neq N$ be a normal subgroup of
$G$. $(G,N)$ is called a \emph{Camina pair} if $xN \subseteq x^G$ for
all $x \in G-N$. $G$ is called a \emph{Camina group} if $(G,\gamma_{2}(G))$ 
is a  Camina pair. The study of such groups was started in \cite{aC78}.

The following lemma is easy to prove:
\begin{lemma}\label{lemma1}
If $(G,N)$ is a Camina pair and $M$ is a normal subgroup of $G$ contained in
$N$ then $(G/M,N/M)$ is a Camina pair.
\end{lemma}

The next theorem follows from \cite{iM81} and \cite{iM86}.
\begin{thm} \label{thm1}
Let $G$ be a finite Camina $p$-group of class $3$ such that 
$[G:  \gamma_{2}(G)] = p^t$, $[\gamma_2(G):\gamma_3(G)] = p^s$ and 
$|\gamma_3(G)| = p^r$. Then
\begin{subequations}
\begin{align}
&\text{$(G,\gamma_{3}(G))$ is  Camina pair, 
$\gamma_{3}(G) = Z(G)$, $t = 2s$, and $s$ is even.}\\
&\text{$s \geq r.$ }\\
&\text{$\gamma_{i}(G)/\gamma_{i+1}(G)$ has exponent $p$, for $i = 1,\;2$.}
\end{align}
\end{subequations}
\end{thm}

\begin{lemma} \label{lemma2}
Let $G$ be a group and $H$ a subgroup of $G$ such that $\gamma_{2}(G) =
\gamma_{2}(H)$. Then $\gamma_{i}(G) = \gamma_{i}(H)$ for all $i \geq 2$.
\end{lemma}
\begin{proof} Since $[G,G]=[H,H] \leq H \leq G$, $H$ is a normal subgroup of 
$G$. 
It follows that the subgroups $\gamma_{i+1}(H) = [\gamma_{i}(H),H]$ are normal
in $G$ for all $i \geq 1$. 
We use induction on $i$ to prove the lemma. 
Since $\gamma_{2}(G) = \gamma_{2}(H)$, assume by induction that 
$\gamma_{i}(G) = \gamma_{i}(H)$, where $i \geq 2$.
We have 
\[[G,\;\gamma_{i-1}(H), H] \leq [\gamma_{i}(G), H] = [\gamma_{i}(H), H] 
= \gamma_{i+1}(H).\]
Since $[H,G] \leq [G,G] = [H,H]$ and $[\gamma_{i}(H), \gamma_{j}(G)] \leq
\gamma_{i+j}(G)$ \cite[Hauptsatz III.2.11]{bH67}, we get
\[[H,G, \gamma_{i-1}(H)] \leq [H, H, \gamma_{i-1}(H)] =
[\gamma_{2}(H),\gamma_{i-1}(H)] \leq \gamma_{i+1}(H).\]
Thus by three subgroup lemma 
\[\gamma_{i+1}(G) = [\gamma_{i}(G),\;G] = [\gamma_{i}(H),\;G] = 
[\gamma_{i-1}(H), H, G] \leq \gamma_{i+1}(H).\]
Since $\gamma_{i+1}(H) \leq \gamma_{i+1}(G)$, it follows that 
$\gamma_{i+1}(G) = \gamma_{i+1}(H)$. This completes the inductive proof of the
lemma. \hfill $\Box$

\end{proof}

\begin{lemma}\label{lemma1a}
Let $G$ be a Camina $p$-group of class $2$, with a
minimal generating set $\{x_1, x_2, \dots, x_t\}$, where $t > 2$.  Then $H =
\gen{x_1, x_2,\dots, x_{t-1}}$ satisfies $\gamma_{2}(G) = \gamma_{2}(H)$.  
In particular, $H$ is a maximal subgroup of $G$.
\end{lemma}
\begin{proof}
Since $G$ is a Camina $p$-group of class $2$, it follows that
$G/\gamma_{2}(G)$ is an elementary abelian $p$-group with
order $p^t$, and that $\gamma_{2}(G) = \Z(G)$ is an elementary abelian 
$p$-group with order $p^s$, for some integer $s > 0$.

First suppose that $s = 1$.  Then $G$ is extra-special.  So $t$ is
even. If $\gamma_{2}(H) \neq \gamma_{2}(G)$, then $\gamma_{2}(H) = 1$, 
since $|\gamma_{2}(G)| = p$ and
$\gamma_{2}(H) \le \gamma_{2}(G)$. In that case $H\Z(G)$ is an abelian 
subgroup of the extra-special group $G$, and so has order 
$\le p^{1 + t/2}$. Since $|H\Z(G)| = p^t$, we get $t \le 1 + t/2$. This 
implies that $t \le 2$, which is a contradiction to our assumption. 
Thus the lemma holds when $s = 1$.

Now suppose that $s > 1$.  If $\gamma_{2}(H) \neq \gamma_{2}(G)$, 
then $\gamma_{2}(H)$ is contained in some maximal subgroup $N$ of 
$\gamma_{2}(G) = \Z(G)$. Then $\bar G = G/N$ is an extra-special group 
with the images $\bar x_1, \bar x_2, \dots, \bar x_t$ of the $x_i$ as 
a minimal generating set.  The subgroup $\gen{\bar x_1, \dots, \bar x_{t-1}}$ 
is the image $\bar H$ of $H$ in $\bar G$.  So $\gamma_{2}(\bar H)$ is the 
image $1$ of $\gamma_{2}(H)$. But $\gamma_{2}(\bar H) = \Z(\bar G) \neq 1$, 
since the lemma holds for $\bar G$. This   contradiction shows that the 
lemma holds for all values of $s$. \hfill $\Box$

\end{proof}

\begin{prop}\label{prop1}
Let $\{x_1, x_2, \dots, x_t\}$ be a minimal generating set
for a Camina $p$-group $G$ of class $3$. Then $H = \gen{x_1, x_2, \dots,
x_{t-1}}$ satisfies $\gamma_{2}(H) = \gamma_2(G)$. In particular, $H$ is a
maximal subgroup of $G$. Moreover, $\gamma_{2}(G) \not\le \Z(H)$.
\end{prop}
\begin{proof}
Set $|G/\gamma_{2}(G)| = p^t$,  $|\gamma_2(G)/\gamma_3(G)| = p^s$ and
$|\gamma_3(G)| = p^r$.  Then it follows from Theorem \ref{thm1} that 
$\gamma_{3}(G)  = \Z(G)$,  $t = 2s$, $s$ is even, and
$s \ge r > 0$. So $t > 2$.  Now it follows from Lemma \ref{lemma1} that 
$G/\gamma_{3}(G)$ is a Camina group. Also the class of $G/\gamma_{3}(G)$ 
is $2$. 
It follows that $G/\gamma_{3}(G) = \gen{x_1\gamma_{3}(G), \dots, 
x_t\gamma_{3}(G)}$ and  $H\gamma_3(G)/\gamma_{3}(G) 
= \gen{x_1 \gamma_{3}(G), \dots, x_{t-1}\gamma_{3}(G)}$.
Thus Lemma \ref{lemma1a} gives that 
$\gamma_{2}(G/\gamma_{3}(G)) = \gamma_{2}(H\gamma_3(G)/\gamma_{3}(G))$. 
This implies that 
\[\gamma_{2}(G)/\gamma_{3}(G) = \gamma_{2}(H)\gamma_{3}(G)/\gamma_{3}(G),\] 
since $\gamma_{3}(G) \le\gamma_{2}(G)$. Thus we have
\begin{equation}\label{e2a}
\gamma_{2}(H)\gamma_3(G) = \gamma_2(G).
\end{equation}

First suppose that $r = 1$, and that $\gamma_{2}(H) < \gamma_2(G)$. 
Since $\gamma_3(G)$ has order $p$, it follows from this and \eqref{e2a} that
$\gamma_2(G)$ is the direct product of its elementary abelian subgroups
$\gamma_{2}(H) \cong \gamma_2(G)/\gamma_3(G)$ and $\Z(G) = \gamma_3(G)$, 
both of which are centralized by $H$. 
So $H$ centralizes $\gamma_2(G)$.  
It then follows from Theorem 1.3(iv) of \cite{MS91} that
$[G:C_{G}(\gamma_2(G))] = p^s$.  Since $\gamma_{2}(G)$ is abelian, 
$H\gamma_2(G) \le C_{G}(\gamma_2(G))$. But $[G:H\gamma_2(G)] = p$, therefore
$p^s = [G:C_{G}(\gamma_2(G))] \le p$. As $s$ is even and $> 0$, this is 
impossible.  Therefore the proposition holds when $r = 1$.

Now assume that $r > 1$, and that $\gamma_2(H) < \gamma_2(G)$. In view
of \eqref{e2a} this implies that $\gamma_2(H) \cap \gamma_3(G) < \gamma_3(G)$. 
So there is some maximal subgroup $N$ of $\gamma_3(G) = \Z(G)$ containing
$\gamma_2(H) \cap \gamma_3(G)$.  Since $N < \gamma_3(G)$, it follows from 
this and Lemma \ref{lemma1} that the factor group $\bar G = G/N$ is a Camina 
$p$-group of class $3$. Since $\gamma_2(\bar G) = \gamma_2(G)/N$
and $\gamma_3(\bar G) = \gamma_3(G)/N$, the images $\bar x_1, \dots,
\bar x_t$ of the $x_i$ form a minimal generating set for $\bar G$, and
$\gen{\bar x_1, \dots, \bar x_{t-1}}$ is the image $\bar H = HN/N$ of $H$ in
$\bar G$.  Notice that $\bar G$ has the same values of $t$ and $s$ as
before, but the new value $1$ for $r$.  Since the proposition holds when
$r = 1$, the image $\gamma_2(\bar H) = \gamma_2(H)N/N$ of $\gamma_2(H)$ is 
equal to $\gamma_2(\bar G) = \gamma_2(G)/N$. Hence 
$\gamma_2(H)N = \gamma_2(G)$. So
$\gamma_3(G) = \gamma_2(H)N \cap \gamma_3(G) = (\gamma_2(H) \cap \gamma_3(G))N 
= N < \gamma_3(G)$, which is impossible. This contradiction proves that
$\gamma_2(H) = \gamma_2(G)$. 

Now it follows from Lemma \ref{lemma2} that $\gamma_3(H) = \gamma_3(G) \neq
1$.  So if $\gamma_2(G) \le \Z(H)$, then $\gamma_2(H) \le \Z(H)$ and therefore
$\gamma_3(G) = \gamma_3(H) = 1$. 
This contradiction proves the proposition.     \hfill $\Box$

\end{proof}

\begin{lemma}\label{lemma2a}
Let $G$ be a Camina $p$-group of class $3$ such that $|\gamma_{3}(G)| \ge
p^2$. Let $H$ be any maximal subgroup of $G$. Then $\Z(H) =\Z(G)$.
\end{lemma}
\begin{proof}Let $H$ be a maximal subgroup of $G$ and $x \in G$ such that 
$G/H = \gen{xH}$. 
Since $Z(G) = \gamma_{3}(G) \le \gamma_2(G) = \Phi(G) \le H$, it follows that 
$\Z(G) \leq \Z(H)$. Suppose that $\Z(G) < \Z(H)$. Then there is 
some non trivial    element $z \in \Z(H)-\Z(G)$. If $z \in
\gamma_2(G)-\gamma_3(G)$, then  $[z,\;G] = \gamma_{3}(G)$ and if $z \in
G-\gamma_2(G)$, then $[z,\;G] =  \gamma_{2}(G)$. In both of the cases it 
follows that 
\[p^2 \le |\gamma_{3}(G)| = |\Z(G)| \leq |[z,\;G]| = |[z, \gen{H,x}]| = 
|[z,\;\gen{x}]| = p,\]
since $x^p \in H$ and $z \in \Z(H)$. This contradiction proves that 
$\Z(G) = \Z(H)$.       \hfill $\Box$

\end{proof}

\section{Almost Camina $p$-groups}

A finite $p$-group $G$ is called an \emph{almost Camina $p$-group} if
$\gamma_{2}(G) = \Phi(G)$ and the following condition (C) holds for every
minimal generating set of $G$:

\parbox{12cm}{(C) \em{If $\{x_1,x_2, \dots, x_d\}$ be a minimal generating 
set for $G$, then  $[x_i,G] = \gamma_{2}(G)$ for all but at most one 
$i = 1, 2, \dots, d$.}}

Notice that every Camina $p$-group is almost Camina. But converse is not true.
For example we take a group $G$ of order $p^4$ and nilpotency class $3$. It is
easy to check that $G$ is an almost Camina $p$-group but it is not a Camina
group. 
\begin{lemma}\label{acaminalemma}
Let $G$ be an almost Camina $p$-group. Let $\{x=x_1, x_2, \dots, x_d\}$ be 
a minimal
generating set for $G$ such that $[x, G] \subset \gamma_{2}(G)$. Then
\begin{subequations}
\begin{align}
&\text{$[y,G] = \gamma_{2}(G)$ for all 
$y \in G-\gen{x}\gamma_{2}(G)$.}\label{acamina1} \\
&\text{$[y,G] = \gamma_{2}(G)$ for all $y \in H-\gamma_{2}(G)$, 
where $H = \gen{x_2, \dots, x_d} \gamma_{2}(G)$.}\label{acamina2} 
\end{align}
\end{subequations}
\end{lemma}
\begin{proof}
Let $V = G/\Phi(G)$. Then $V$ is a vector space of dimension $d$ over $\mathbb
F_{p}$ and $\{v_1, \dots, v_d\}$ is a basis for $V$, where 
$v_i = x_i\Phi(G)$, $1 \le i \le d$.  Let $W$ be the subspace of $V$ spanned 
by $v_1$. Now suppose that there is an element $y$ in $G-\gen{x}\gamma_{2}(G)$ 
such that $[y,G] \subset \gamma_{2}(G)$. Since $y \in G-\gen{x}\gamma_{2}(G)$, 
$w = y\Phi(G) \in V-W$. Thus the set $\{v_1, w\}$ is linearly independent. 
Therefore we can extend this set to a basis $\{v_1, w = w_2, \dots, w_d\}$ 
of $V$. Hence both $x$ and $y$ can
be included in some minimal generating set of $G$.
This contradicts the fact that $G$ is an almost Camina group, since both
$[x,G]$ as well as $[y,G]$ are proper subsets of $\gamma_{2}(G)$. This proves
\eqref{acamina1}.

Now consider $H$ as in \eqref{acamina2}. Let $y \in H -\gamma_{2}(G)$.
We claim that $y \in G-\gen{x}\gamma_{2}(G)$. Suppose that 
$y \in \gen{x}\gamma_{2}(G)$. Then $y = x^i u$ for some $u \in \gamma_2(G)$
and some integer $i$ such that $1 \le i \le p-1$, since 
$y \not\in  \gamma_2(G)$. This implies that $x^i = yu^{-1} \in H$, since
$u \in \gamma_2(G) \le H$ and $y \in H$. Hence $x \in H$ and therefore $H=G$,
a contradiction. Thus our claim holds. Now  \eqref{acamina2} follows from
\eqref{acamina1}.            \hfill $\Box$

\end{proof}

\begin{prop}\label{prop1a}
Let $G$ be an almost Camina $p$-group  such that $|G| = p^{2m}$ 
and $|\gamma_{2}(G)| = p^m$, where $m \ge 3$. 
Then there exists a maximal subgroup $H$ of $G$ such that 
$\gamma_{2}(G) = \gamma_{2}(H)$. Moreover, if the nilpotency 
class of $G$ is at least $3$, then $\gamma_{2}(G) \not\le \Z(H)$ for any such 
maximal  subgroup $H$.
\end{prop}
\begin{proof} 
If $G$ is a Camina $p$-group, then the proposition holds from Lemma
\ref{lemma1a} and Proposition \ref{prop1}. So we assume that $G$ is not a
Camina group. We complete the proof of this proposition in three steps.
\setcounter{step}{0}
\begin{step}\label{step1a}
Let $x' \in G-\gamma_{2}(G)$ such that  
$\C_G(x')\gamma_{2}(G) \neq \gen{x'}\gamma_{2}(G)$. 
Then $\gamma_2(G) = \gamma_2(H)$ for some maximal subgroup $H$ of $G$.
\end{step}
\begin{proof}
Since $\C_G(x')\gamma_2(G) \neq \gen{x'}\gamma_2(G)$, there exists an element
$y' \in G- \gen{x'}\gamma_{2}(G)$ such that $x'y' = y'x'$. 
Then $M=\gen{x', y'}$ is an abelian subgroup of $G$ and the order of
$M\gamma_{2}(G)/\gamma_{2}(G)$ is $p^2$. 
Choose $x \in M-\gamma_{2}(G)$ such that $|[x,G]|$ is minimum, i.e., 
$|[x,G]| \le |[z,G]|$ for all $z \in M-\gamma_{2}(G)$. 
Let $\{x = x_1, x_2, \dots, x_{m}\}$ be a minimal generating set
for $G$. Let $H = \gen{x_2, \dots, x_{m}}\gamma_{2}(G)$. 
Since $x^p \in \gamma_{2}(G)$, the order of $\gen{x}\gamma_{2}(G)/\gamma_2(G)$ 
is $p$. Thus $\gen{x}\gamma_{2}(G) < M \gamma_2(G)$. 
Let $y \in M-\gen{x}\gamma_{2}(G)$. Then $xy=yx$.
Now either $[x,G] = \gamma_{2}(G)$ or $[x,G] \subset \gamma_{2}(G)$. 
If $[x,G] = \gamma_{2}(G)$, then by  the minimality of $|[x,G]|$
it follows that $[y,G] = \gamma_{2}(G)$. And, if 
$[x,G] \subset \gamma_{2}(G)$, then it follows  from Lemma 
\ref{acaminalemma} that  $[y, G] = \gamma_{2}(G)$. Thus in both of the two
cases we have $[y, G] = \gamma_{2}(G)$.
Let $u$ be an arbitrary element of $\gamma_2(G)$. Then $u = [y, g]$ for 
some $g \in G$. 
Now $g$ can be written as $g = hx^i$ for some $h \in H$ and some integer $i$. 
Therefore
\[u = [y, g] = [y,  hx^i] = [y, x^i][y, h]^{x^i} = [y, h]^{x^i} \in
\gamma_{2}(H),\] as $y \in H$ and  $\gamma_{2}(H)$ is a normal subgroup of 
$G$.  This implies that 
$\gamma_{2}(G) \leq \gamma_{2}(H)$. Since the other way inclusion is obvious, 
we have  $\gamma_2 (G) = \gamma_2 (H)$. 
\end{proof}

\begin{step}\label{step2a}
Let $x \in G-\gamma_{2}(G)$ such that $|[x,G]| \le p^{m-2}$.
Then $\gamma_2(G) = \gamma_2(H)$ for some maximal subgroup $H$ of $G$.
\end{step}
\begin{proof}
Since $|[x,G]| \le p^{m-2}$, therefore $|\C_{G}(x)| \ge p^{m+2}$. But 
$|\gen{x}\gamma_{2}(G)| = p^{m+1}$. Therefore 
$\C_G(x)\gamma_{2}(G) \neq \gen{x}\gamma_{2}(G)$. Now Step \ref{step2a} follows
from Step \ref{step1a}.
\end{proof}

\begin{step}\label{step3a}
Let $\C_G(u)\gamma_{2}(G) = \gen{u}\gamma_{2}(G)$ for all $u \in
G-\gamma_{2}(G)$. Then it follows from Step \ref{step2a} that $|[u,G]| \geq
p^{m-1}$ for all $u \in G-\gamma_{2}(G)$. Let $x \in G-\gamma_{2}(G)$ 
such that $|[x,G]| = |\gamma_{2}(G)|/p = p^m/p =p^{m-1}$. Let 
$\{x = x_1, x_2, \dots, x_{m}\}$ be a minimal generating set for $G$.
Then $H = \gen{x_2, \dots, x_{m}}\gamma_{2}(G)$ is a maximal subgroup of $G$
such that $\gamma_2 (G) = \gamma_2 (H)$.
\end{step}
\begin{proof}
Set $N = \gen{x}\gamma_{2}(G)$. Notice that $|N| = p^{m+1}$. Then it 
follows from Lemma \ref{acaminalemma} that for any 
$y \in G-N$, $[y,G] = \gamma_{2}(G)$. Since 
$\C_G(u)\gamma_{2}(G) = \gen{u}\gamma_{2}(G)$ for all 
$u \in G-\gamma_{2}(G)$, we have $\C_G(u) \le \gen{u}\gamma_{2}(G)$. In
particular, $\C_{G}(x) \le \gen{x} \gamma_{2}(G)  = N$. Since
$|x^G| = |[x,G]| = p^{m-1}$, therefore $|\C_{G}(x)| = |G|/|x^G| = p^{m+1}$. 
Thus
$\C_{G}(x) = N$, since $|N| = p^{m+1}$. This implies that $\gamma_{2}(G) \le
\C_{G}(x)$. Now both $[x,G]$ as well as $\gamma_{2}(H)$ are subsets of 
$\gamma_{2}(G)$ containing $1$, therefore, either 
$[x,G] \cap \gamma_{2}(H) \neq 1$ or $[x,G] \cap \gamma_{2}(H) =1$.

First assume that $[x,G] \cap \gamma_{2}(H) \neq 1$. Let $1 \neq [x,g] \in 
[x,G] \cap \gamma_{2}(H)$. Then $[x,g] = [x,hx^k]$ for some $h \in
H-\gamma_{2}(G)$ and some integer $k$ such that $1 \le k \le p$. Therefore 
\[[x,g] = [x,hx^k] = [x,x^k][x,h]^{x^k} = [x,h],\]
since $\gamma_{2}(G) \le \C_{G}(x)$. Since $h \in H-\gamma_{2}(G)$, it again
follows from Lemma \ref{acaminalemma} that $[h,G] = \gamma_{2}(G)$. 
Let $[h,g']$ be an arbitrary element of $\gamma_{2}(G)$. Then
$[h,g'] = [h, h'x^i]$ for some $h' \in H$ and some integer $i$ such that 
$1 \le i \le p$. Thus we have that 
\[[h,g'] = [h,h'x^i] = [h,x^i][h,h']^{x^i} 
= [h,x]^i[h,h'] \in \gamma_{2}(H),\]
 since $\gamma_{2}(G) \le \C_{G}(x)$
and $[h,x] \in \gamma_{2}(H)$. This implies that $\gamma_{2}(G) \le
\gamma_{2}(H)$. Since $\gamma_{2}(H) \le \gamma_{2}(G)$, we have 
$\gamma_{2}(H) = \gamma_{2}(G)$.

Now assume that $[x,G] \cap \gamma_{2}(H) = 1$. This implies that
$\gamma_{2}(H) < \gamma_{2}(G)$. Let $h \in H-\gamma_{2}(G)$. Then
$\C_{G}(h) \le \gen{h}\gamma_{2}(G) \le H$. Thus $\C_{H}(h) = \C_{G}(h) \cap H
= \C_{G}(h)$. Therefore 
\[|[h,H]| = |h^H| = |H|/|\C_{H}(h)| = |H|/|\C_{G}(h)| = p^{2m-1}/p^m 
= p^{m-1},\]
since 
\[|\C_{G}(h)| = |G|/|h^G| = |G|/|[h,G]| = |G|/|\gamma_{2}(G)| = p^{2m}/p^m
= p^m.\]
This proves that $[h,H] = \gamma_{2}(H)$, since $[h,H] \le \gamma_{2}(H)$ and 
$\gamma_{2}(H) < \gamma_{2}(G)$. Thus $|\gamma_{2}(H)| = p^{m-1}$.

Now consider the factor group $\bar{G} = G/\gamma_{2}(H)$. Since
$\gamma_{2}(G) \not\le \gamma_{2}(H)$, $\bar{G}$ is a non-abelian group of 
order $p^{m+1}$. We have $\gamma_{2}(\bar{G}) = \overline{\gamma_{2}(G)}$, 
since  $\gamma_{2}(H) \le \gamma_{2}(G)$. 
Now $|\gamma_{2}(\bar{G})| = |\gamma_{2}(G)|/
|\gamma_{2}(H)| = p$. Thus $\gamma_{2}(\bar{G}) \leq \Z(\bar{G})$. Therefore
the nilpotency class of $\bar{G}$ is $2$. Now $[\bar{g}, \bar{G}] =
\gamma_{2}(\bar{G})$ for all $\bar{g} \in \bar{G}- \gamma_{2}(\bar{G})$. 
For, if there is some element $\bar{g} \in \bar{G}- \gamma_{2}(\bar{G})$ such
that $[\bar{g}, \bar{G}] < \gamma_{2}(\bar{G})$, then $|[\bar{g}, \bar{G}]| <
p$. Thus $|[\bar{g}, \bar{G}]| = 1$ and $[g,G]\gamma_{2}(H) =
\gamma_{2}(H)$. This implies that $[g,G] \subseteq \gamma_{2}(H)$. If $g \in
G-N$, then $\gamma_{2}(G) = [g,G] \le \gamma_{2}(H)$, which 
is a contradiction to the fact that $\gamma_{2}(H) < \gamma_{2}(G)$. So let
$g \in N-\gamma_{2}(G)$. Then $g = vx^i$ for some $v \in \gamma_{2}(G)$ and
some integer $i$ such that $1 \le i \leq p-1$. Since $[g,G] \subseteq
\gamma_{2}(H)$, therefore $[g,H] \subseteq \gamma_{2}(H)$. 
Let $h$ be any element of $H-\gamma_{2}(G)$. Then $[g,h] \in \gamma_2(H)$ and 
\[[g,h] = [vx^i,h] = [v,h]^{x^i}[x^i,h] = [v,h][x,h]^i \in \gamma_{2}(H).\]
Now $v \in \gamma_{2}(G) \le H$, therefore $[v,h] \in \gamma_{2}(H)$. This
implies that $[x,h]^i \in \gamma_{2}(H)$ and therefore $[x,h] \in
\gamma_{2}(H)$. Then $[x,h]$ is either $1$ or not. If $[x,h] = 1$, then
$h \in H-\gamma_{2}(G)$ centralizes $x$, so that $C_{G}(x)\gamma_2(G) \neq 
\gen{x}\gamma_{2}(G)$, contrary to the hypotheses of this step. If $[x,h] \neq
1$, then $1 \neq [x,h] \in [x,G] \cap \gamma_{2}(H)$ contradicts the
assumption that $[x,G] \cap \gamma_2(H) = 1$.
Thus $\bar{G}$ is a Camina group of class $2$. It then follows from
Corollary 2.4 of \cite{iM81} that $\bar{G}$ is special, therefore $\bar{G}$ 
is extra special. Thus $m$ must be even. Now $\bar{H}$ is an abelian
normal subgroup of $\bar{G}$. Thus $|\bar{H}| \le p^{1+\frac{m}{2}}$. But
$|\bar{H}| = p^m$, therefore $p^m \le p^{1+\frac{m}{2}}$. This implies that
$m \le 1+\frac{m}{2}$. Thus $m \le 2$, which is a contradiction to our
assumption that $m \ge 3$. Hence the latter case can not occur. This completes
the proof of Step \ref{step3a}.  
\end{proof}
Thus the existence of a maximal subgroup $H$ of $G$ such that $\gamma_{2}(G)
=\gamma_{2}(H)$ follows from Step \ref{step1a} -- Step \ref{step3a}. 
Let $H$ be any maximal subgroup of $G$ such that $\gamma_{2}(G)
=\gamma_{2}(H)$. Now it follows from Lemma \ref{lemma2} that the nilpotency 
class of $H$ is equal to the nilpotency class of $G$. 
If $\gamma_2(G) \le \Z(H)$,  then $\gamma_{2}(H) = \gamma_2(G) \le \Z(H)$. 
Hence the class of $H$ is $2$ and therefore the class of $G$ is also $2$. 
Thus it follows that $\gamma_{2}(G) \not\le \Z(H)$ whenever the nilpotency 
class of $G$ is $\ge 3$. 
This completes the proof of the proposition. \hfill $\Box$

\end{proof}

\section{Proof of the main theorem}

Let $G$ be a finite $p$-group of order $p^n$.
Let $\{x_{1}, \cdots, x_{d}\}$ be any minimal generating set for $G$. 
Let $\alpha \in \Aut_{c}(G)$. Since $\alpha(x_{i}) \in x_{i}^G$ for 
$1 \leq i \leq d$, there are at the most $|x_{i}^{G}|$ choices for the image
of $x_i$ under $\alpha$. Thus it follows that
\begin{equation}
\label{bineq} |\Aut_{c}(G)| \leq \prod_{i=1}^{d} |x_{i}^G|.
\end{equation}

Let  $|\gamma_{2}(G)| = p^m$. Let $\Phi(G)$ denotes the Frattini subgroup 
of $G$. Since $\gamma_{2}(G)$ is contained in $\Phi(G)$, by the  Burnside 
basis theorem we have $d \leq n-m$. 
Notice that $|x_{i}^G| \leq |\gamma_{2}(G)| = p^m$ for all $i = 1, 2, 
\cdots, d$. So from \eqref{bineq} we get
\begin{equation}
\label{lbineq} |\Aut_{c}(G)| \leq p^{md} \leq (p^{m})^{n-m} = p^{m(n-m)}.
\end{equation}

\begin{thm}\label{thm2}
Let $G$ be a finite $p$-group. 
If equality holds in \eqref{lbineq}, then $G$ is either an abelian $p$-group,
or a non-abelian  Camina special $p$-group. 
\end{thm}
\begin{proof} 
 Let $G$ be a finite $p$-group such that  equality holds in \eqref{lbineq}. 
If $G$ is abelian, then we are done. So assume that $G$ is non-abelian.
Equality in \eqref{lbineq} implies that $d = n-m$ and that
\begin{equation}
\label{lbeq} |\Aut_{c}(G)| = p^{md} = p^{m(n-m)}.
\end{equation}
Since $G/\gamma_2(G)$ is the direct product of $d$ non-trivial cyclic
$p$-groups, yet has order $p^{n-m} = p^d$, it must be elementary abelian. 
It then follows that any arbitrary  element $x \in G - \gamma_{2}(G)$ is a 
part of a minimal generating set $\{x = x_1, x_2, \dots, x_{n-m}\}$ for $G$. 
Now $x[x,G] = x^G$ implies $|[x,G]| = |x^G|$.
If $[x,G] \subset [G,G]$, then $|x^G| = |[x,G]| < |[G,G]| = p^m$ and
\eqref{bineq} gives $|\Aut_c(G)| < p^{m(n - m)}$, which contradicts 
\eqref{lbeq}. 
Thus it follows that $[x,G] = [G, G]$ for all $x \in G - [G,G]$. 
This proves that $G$ is a Camina group. That the
nilpotency class of $G$ is $\leq 3$ follows from Main Theorem of \cite{DS96}. 
  
 It follows from  \eqref{lbeq} that  given any minimal set of generators
$\{x_1, x_2, \dots,$ $ x_{n-m}\}$ for $G$, and any elements $y_1, y_2, \dots,
y_{n-m}$ $\in \gamma_{2}(G)$ (need not be distinct), there is some 
automorphism  $\alpha \in \Aut_c(G)$ such that $\alpha(x_i) = x_iy_i$ for 
$i= 1, 2,\dots,n-m$.  In particular, we can choose $\alpha$ such that
\[\alpha(x_{i}) = x_{i},\;1 \leq i \leq n-m-1 \;\;\text{and} \;\; 
\alpha(x_{n-m}) = x_{n-m}y,\]
where $y$ is  an arbitrary element of $\gamma_{2}(G)$.

     Let $G$ have class 3. Let $G = \gen{x_1,x_2, \dots, x_{n-m}}$. Then
it follows from Proposition \ref{prop1} that 
 $H = \gen{x_1, \dots, x_{n-m-1}}$ is a maximal subgroup of $G$ such that 
$\gamma_{2}(H) = \gamma_{2}(G)$ and $\gamma_{2}(G) \not \leq \Z(H)$. 
We fix an element $y \in \gamma_{2}(G) - \Z(H)$. 
Now we can choose an automorphism $\alpha \in \Aut_c(G)$ such that 
$\alpha(x_{i}) = x_{i}$ for $1 \leq i \leq n-m-1$ and $\alpha(x_{n-m}) 
= x_{n-m}y$.
Then $\alpha$ centralizes $H$, as well
as $G/H = \gen{x_{n-m}H}$.  But any automorphism of $G$ centralizing both
$G/H$ and $H$ must send the generator $x_{n-m}$ of $G$ modulo $H$ 
to $x_{n-m}z$, for some element $z \in \Z(H)$ \cite[Satz I.17.1]{bH67}. 
Since $\alpha(x_{n-m}) = x_{n-m}y$ and $y \notin \Z(H)$,
this is impossible.  Thus the nilpotency class of $G$ can not be $3$. 
Since $G$ is non-abelian,
this proves that the nilpotency class of $G$ is $2$. That $G$ is special 
follows  from Corollary 2.4 of \cite{iM81}. \hfill $\Box$

\end{proof}

Let $G$ be a finite $p$-group of class $2$. Let $\phi \in
\Aut_c(G)$. Then the map $g \mapsto g^{-1}\phi(g)$ is a homomorphism of
$G$ into $\gamma_2(G)$. This homomorphism sends $Z(G)$ to $1$. So it
induces a homomorphism $f_{\phi} \colon G/Z(G) \to \gamma_2(G)$, sending
$gZ(G)$ to $g^{-1}\phi(g)$, for any $g \in G$.  It is easily seen that
the map $\phi \mapsto f_{\phi}$ is a monomorphism of the group
$\Aut_c(G)$ into $\Hom(G/Z(G), \gamma_2(G))$.

Any $\phi \in \Aut_c(G)$ sends any $g \in G$ to some $\phi(g)
\in g^G$. Then $f_{\phi}(gZ(G)) = g^{-1}\phi(g)$ lies in $g^{-1}g^G =
[g,G]$.  Denote
\[  \{ \, f \in \Hom(G/Z(G), \gamma_2(G)) \mid f(gZ(G)) \in [g,G], \text{ for
all $g \in G$}\,\} \]
by $\Hom_c(G/Z(G), \gamma_2(G))$. Then $f_{\phi} \in \Hom_c(G/Z(G),
\gamma_2(G))$ for all $\phi \in \Aut_c(G)$. On the other hand, if $f \in
\Hom_c(G/Z(G), \gamma_2(G))$, then the map sending any $g \in G$ to
$gf(gZ(G))$ is an automorphism $\phi \in \Aut_c(G)$ such that $f_{\phi}
= f$. Thus we have

\begin{prop}\label{prop2}
 Let $G$ be a finite $p$-group of class 2. Then the
above map $\phi \mapsto f_{\phi}$ is an isomorphism of the group
$\Aut_c(G)$ onto $\Hom_c(G/Z(G), \gamma_2(G))$.
\end{prop}

The next lemma follows from \cite[page 335]{jR95}.
\begin{lemma}\label{hlemma}
Let $H$  and $K$ be two elementary abelian groups of order $p^r$ and $p^s$
respectively. Let $H = \times_{i=1}^r H_i$ and $K = \times_{j=1}^s K_j$,
where $H_i$ and $K_j$ are cyclic groups of order $p$, $1 \le i \le r$, $1 \le
j \le s$. Then
\[\Hom(H, K) \cong \times_{i=1, j=1}^{r, s} \Hom(H_i,K_j).\]
In particular, $\Hom(H, K)$ is an elementary abelian group of order
$p^{rs}$.
\end{lemma}

\begin{thm} \label{thm3}
Let $G$ be a finite $p$-group. Equality holds in \eqref{lbineq} if and only
if $G$ is either an abelian $p$-group, or a non-abelian  Camina special 
$p$-group.
\end{thm}
\begin{proof}  Let $G$ be finite $p$-group such that $|G| = p^n$ and 
$|\gamma_{2}(G)| = p^m$. 
Let $G$ be abelian. 
Then $|\gamma_{2}(G)| = 1$ implies $m = 0$.
We have $|\Aut_{c}(G)| = 1 = p^{m(n-m)}$. Thus equality holds in
\eqref{lbineq}.
Now let  $G$ be a non-abelian Camina special $p$-group. Then 
 $[x,G] = \gamma_{2}(G)$ for all non-central elements $x \in G$ and 
$[x,G][y,G] = [xy,G]$ for all $x,y \in G - \Z(G)$ such that $xy \not\in \Z(G)$.
Pick any $f \in \Hom(G/\Z(G), \gamma_{2}(G))$. Then 
\[f(\overline{x}) \in \gamma_{2}(G) = [x,G] = x^{-1}x^G,\]
where $\overline{x} = x\Z(G)$. 
Thus $f \in \Hom_{c}(G/\Z(G), \gamma_{2}(G))$ and therefore 
\[\Hom(G/\Z(G), \gamma_{2}(G)) \leq \Hom_{c}(G/\Z(G), \gamma_{2}(G)).\]
Since $\Hom_{c}(G/\Z(G), \gamma_{2}(G)) \leq \Hom(G/\Z(G), \gamma_{2}(G))$, we
have 
\[\Hom_{c}(G/\Z(G), \gamma_{2}(G)) = \Hom(G/\Z(G), \gamma_{2}(G)).\]
Since $G$ is a special $p$-group, it follows that $\Z(G)$ and $G/\gamma_2(G)$
are elementary abelian groups of order $p^m$ and $p^{n-m}$ respectively.
Thus from Proposition \ref{prop2} and Lemma \ref{hlemma}, we get 
\[|\Aut_{c}(G)| = |\Hom(G/\Z(G), \gamma_{2}(G))| = p^{m(n-m)}.\] 

The converse follows from Theorem \ref{thm2}. \hfill $\Box$

\end{proof}

Now we prove Theorem A.
\begin{thm} \label{thm4}
Let $G$ be a non-trivial $p$-group having order $p^{n}$. Then
\[|\Aut_{c}(G)| \leq
   \begin{cases}
     p^{\frac{(n^{2}-4)}{4}},  &\text{if $n$ is even;}\\
     p^{\frac{(n^{2}-1)}{4}}, &\text{if $n$ is odd.}
   \end{cases}
\]
\end{thm}
\begin{proof}  If $G$ is abelian, then the theorem holds trivially. So let $G$
be non-abelian and $|\gamma_{2}(G)| = p^{m}$. Notice that 
$|x^G| \leq |\gamma_{2}(G)| = p^{m}$
 for all $x \in G$. Let $|\Phi(G)| = p^{t}$. Since $\gamma_{2}(G) \subseteq
 \Phi(G)$, $m \leq t$. By the Basis Theorem of Burnside it follows that 
from any
 generating set for $G$ one can choose  $n-t$ elements such that these $n-t$
 elements generate $G$. $n-t$ is maximum if $t = m$. 
Since $|\gamma_{2}(G)| = p^m$, we have  $1 \leq m \leq n-2$. 
Thus all possible values of $m(n-m)$ are
\[\{n-1,\; 2(n-2),\; 3(n-3),\;\cdots , \;n^{2}/4\}\]
if $n$ is even and
\[\{n-1,\; 2(n-2),\; 3(n-3),\;\cdots , \;(n-1)(n+1)/4\}\]
if $n$ is odd. Notice that the maximum value of $m(n-m)$ is $n^{2}/4$
when $n$ is even, and $(n-1)(n+1)/4$ when $n$ is odd. Putting these 
values in formula  \eqref{lbineq} we get
\[|\Aut_{c}(G)| \leq
   \begin{cases}
     p^{\frac{n^{2}}{4}},  &\text{if $n$ is even;}\\
     p^{\frac{(n^{2}-1)}{4}}, &\text{if $n$ is odd.}
   \end{cases}
\]
Thus in the case when $n$ is odd we are done.

Now assume that $n$ is even. Let $|\Aut_{c}(G)| = p^{n^{2}/4}$. This is
possible only when $m = n/2$. So assume that $m = n/2$. 
Thus $|\Aut_{c}(G)| =  p^{m(n-m)}$ and equality holds in \ref{lbineq}. 
Now by Theorem \ref{thm2} it follows  that $G$ is a Camina special
$p$-group. 
It then follows from Theorem 3.2 of \cite{iM81} that $n-m$  is even and 
$n-m \geq 2m$. 
This implies that 
\begin{equation}
\label{goodineq} m \leq n/3, 
\end{equation}
which contradicts our assumption that $m = n/2$.
Thus there exists no finite $p$-group $G$ such that 
$|\Aut_{c}(G)| = p^{n^{2}/4}$.
Therefore we have $|\Aut_{c}(G)| < p^{n^{2}/4}$. Thus
\[|\Aut_{c}(G)| \leq p^{\frac{n^{2}}{4} - 1} = p^{\frac{(n^2 - 4)}{4}}.\]
This, along with the case when $n$ is odd, proves the theorem. \hfill $\Box$

\end{proof}

\begin{lemma}\label{ecd1}
Let $G$ be a non-abelian group of order $p^{2m}$ such that 
$|\gamma_{2}(G)| = p^m$ and  equality holds in \eqref{gbineq}. 
Then $G$ is an almost Camina group of nilpotency class $\ge 3$, which is not
a Camina group. Moreover, $|[x,G]| \ge p^{m-1}$ for all $x \in G-\gamma_2(G)$.
\end{lemma}
\begin{proof}
Since equality holds in \eqref{gbineq}, $\gamma_{2}(G) = \Phi(G)$. For, if
$\gamma_{2}(G) < \Phi(G)$, then it follows from \eqref{bineq} that
$|\Aut_{c}(G)| \le p^{m(m-1)} < p^{m^2-1}$ for all $m \ge 2$. Since $G$ is
non-abelian, $m \ge 2$. This contradicts our hypothesis that equality 
holds in \eqref{gbineq}. If the nilpotency class of $G$ is $2$, then
$\gamma_{2}(G) \le \Z(G)$. Thus for any element $x \in G-\gamma_{2}(G)$, 
$|\C_{G}(x)| \geq p^{m+1}$. Therefore $|x^G| \leq p^{m-1}$ for all
$x \in G-\gamma_{2}(G)$. Then again it follows from \eqref{bineq} that
$|\Aut_{c}(G)| \le p^{m(m-1)} < p^{m^2-1}$. This again gives a contradiction
to our supposition. Thus nilpotency class of $G$ is $\ge 3$.

Now we claim that $G$ is not a Camina group. Assume the contrary, i.e., 
$G$ is a Camina group of class $3$. Let 
$[G:  \gamma_{2}(G)] = p^t$, $[\gamma_2(G):\gamma_3(G)] = p^s$ and 
$|\gamma_3(G)| = p^r$. Then it follows from Theorem \ref{thm1} that 
$\gamma_{3}(G) = Z(G)$, $t=2s$, $s$ is even, and $s \geq r$. Then
\[p^{2m} = |G| = p^{t+s+r}.\]
Since $t$ and $s$ are even, it follows that $r$ is even. Thus $r \ge 2$.
Let $\{x_1, \dots, x_m\}$ be a minimal generating set for $G$. Then it
follows from Proposition \ref{prop1} that $H =\gen{x_1, \dots, x_{m-1}}$
is a maximal subgroup of $G$ such that $\gamma_{2}(G) = \gamma_{2}(H)$. Then
$G/H = \gen{x_mH}$. Let $A$ be the subgroup of $\Aut_{c}(G)$ consisting of
all $\alpha \in \Aut_{c}(G)$ which centralize both $H$ and $G/H$. 
If $\alpha \in A$, then it
follows from \cite[Satz I.17.1]{bH67} that $\alpha(x_m) = x_mz$ for some
$z \in \Z(H)$.  
Since equality holds in \eqref{gbineq}, it follows that the orbit
of $x_m$ under the action of $A$ must have length at least $p^{m-1}$. Thus
$\Z(H) \ge p^{m-1}$. But from Lemma \ref{lemma2a} we have that 
$\Z(G) = \Z(H)$, since $r \ge 2$. So if $|\Z(H)| = p^{m-1}$, then 
\[[\gamma_{2}(G):\gamma_{3}(G)] = [\gamma_{2}(G):\Z(G)] 
= [\gamma_{2}(G): \Z(H)] = p.\]
This implies that $s = 1$, which is a contradiction to the fact that $s$ is
even. And if $|\Z(H)|=p^m$, then $\gamma_{2}(G) = \Z(G)$. Since the class of 
$G$ is $3$, it is impossible. Hence our claim is true.

So we assume that $G$ is not a Camina group. Then there must exist some
minimal generating set $\{x_1, \dots, x_m\}$ for $G$ such that
$[x_i,G] \subset \gamma_{2}(G)$  for at least one $i$, $1 \le i \le m$.
Now assume that $\{x_1, \dots, x_m\}$ is any such minimal generating 
set for $G$. Then we claim that
$[x_i,G] \subset \gamma_{2}(G)$ for at most one $i$, $1 \le i \le m$ and for
this $i$, $|[x_i,G]| = p^{m-1}$. First assume that there are more than one $i$,
$1 \le i \le m$ such that $[x_i,G] \subset \gamma_{2}(G)$. 
So $|[x_i,G]| < p^m$ for more than one $i$. Then it follows from the 
estimates in \eqref{bineq} that $|\Aut_{c}(G)| \le p^{m(m-2)}p^{m-1}p^{m-1}$.  
From equality in \eqref{gbineq} we have that $|\Aut_{c}(G)| = p^{m^2-1}$. 
Thus 
\[p^{m^2-1} = |\Aut_{c}(G)| \le p^{m(m-2)+m-1+m-1}.\]
This implies that $m^2 - 1 \le m^2 - 2m + m-1 + m-1$ or equivalently
$1 \ge 2$, which is absurd. Next assume that for some $i$, 
$|[x_i,G]| \le p^{m-2}$. Then again from the estimates in \eqref{bineq} we
have that $|\Aut_{c}(G)| \le p^{m(m-1)}p^{m-2}$. Thus we have
\[p^{m^2-1} = |\Aut_{c}(G)| \le p^{m(m-1)+m-2}.\] 
This implies that $m^2 -1 \le m^2 -m +m -2$ or equivalently $1 \ge 2$, which is
again absurd. Thus our claim holds.
Now it follows from the definition of almost Camina $p$-groups that $G$ is 
an almost Camina group. \hfill $\Box$

\end{proof}

Let $X$ be a finite group and $\bar{X} = X/\Z(X)$. 
Then commutation in $X$ gives a well defined map
$a_{X} : \bar{X} \times \bar{X} \mapsto \gamma_{2}(X)$ such that
$a_{X}(x\Z(X), y\Z(X)) = [x,y]$ for $(x,y) \in X \times X$.
Two finite groups $G$ and $H$ are called \emph{isoclinic} if 
there exists an  isomorphism $\phi$ of the factor group
$\bar G = G/\Z(G)$ onto $\bar{H} = H/\Z(H)$, and an isomorphism $\theta$ of
the subgroup $\gamma_{2}(G)$ onto  $\gamma_{2}(H)$
such that the following diagram is commutative
\[
 \begin{CD}
   \bar G \times \bar G  @>a_G>> \gamma_{2}(G)\\
   @V{\phi\times\phi}VV        @VV{\theta}V\\
   \bar H \times \bar H @>a_H>> \gamma_{2}(H).
  \end{CD}
\]
The resulting pair $(\phi, \theta)$ is called an \emph{isoclinism} of $G$ 
onto $H$. Notice that isoclinism is an equivalence relation among finite 
groups.

For the statement of the next proposition
we need the following group of order $p^6$, which is the group 
$\phi_{21}(1^6)$ in the isoclinism family ($21$) of \cite{rJ80}:
\begin{eqnarray} 
R &=& \langle \alpha, \alpha_1, \alpha_2, \beta, \beta_1, \beta_2| 
[\alpha_1, \alpha_2] = \beta, [\beta, \alpha_i]=\beta_i, 
[\alpha, \alpha_1]= \beta_2,\label{r}\\
& &\;\; [\alpha, \alpha_2] = \beta_1^{\nu},
\alpha^p=\beta^p=\beta_i^p = 1, \alpha_1^p=\beta_1^{(^p_3)},
\alpha_2^p=\beta_1^{-(^p_3)},\nonumber\\
& & \;\;i=1,2 \rangle,\nonumber
\end{eqnarray}
where $\nu$ is the smallest positive integer which is a non-quadratic residue
{$\mod p$} and $\beta_1$ and $\beta_2$ are central elements. 
Notice that the subgroup $\gamma_2(R) = \gen{\beta, \beta_1, \beta_2}$ is 
elementary abelian of order $p^3$ and the subgroup $\gamma_3(R) = \Z(R) = 
\gen{\beta_1, \beta_2}$ is elementary abelian of order $p^2$. The quotient
group  $R/\gamma_2(R)$ is elementary abelian of order $p^3$ with
$\{\bar{\alpha}, \bar{\alpha_1}, \bar{\alpha_2}\}$ as a minimal basis, where
$\bar{\alpha} = \alpha \gamma_2(R)$ and $\bar{\alpha_i} = \alpha_i
\gamma_2(R)$, $i = 1, 2$ and the quotient group $\gamma_2(R)/\gamma_3(R) =
\gen{\beta \gamma_3(R)}$ is cyclic of order $p$.

\begin{prop}\label{prop2a}
Let $G$ be an almost Camina group of order $p^{2m}$, $m \ge 2$ such that 
$|\gamma_{2}(G)| = p^m$. Let equality hold in \eqref{gbineq}. 
Then one of the following holds:
\begin{subequations}
\begin{align}
&\text{$G$ is a group of order $p^4$ and nilpotency class $3$;}\\
&\text{$G$ is isoclinic to $R$.}
\end{align}
\end{subequations}
\end{prop}
\begin{proof}
Since $|\gamma_{2}(G)|=p^m$ and $m \geq 2$, $G$ must be non-abelian. 
Thus it follows from Lemma \ref{ecd1} that the nilpotency class of $G$ 
is at least $3$ and $G$ is not a Camina group.
First suppose that $m = 2$. Then obviously $|G| = p^4$ and the nilpotency class
of $G$ is $3$.

Now suppose that $m \geq 3$.  
Since $G$ is almost Camina, but not
Camina, there exists  some element $z \in G - \gamma_2(G)$ such that $[z,G]
\subset \gamma_2(G)$.  Then $[y, G] = \gamma_2(G)$ for all $y \in G -
\gen{z}\gamma_2(G)$ by Lemma \ref{acaminalemma}.  
From Lemma \ref{ecd1} we have that $|[x,G]| \geq p^{m-1}$ for all 
$x \in G-\gamma_{2}(G)$. Thus we have $|[z,G]| = p^{m-1}$. 
It also follows from here that $\Z(G) < \gamma_2(G)$.
Since equality holds in \eqref{gbineq}, given any minimal
generating set for $G$ of the form $\{x_1 = z, x_2, \dots, x_m\}$, any
$y_1 \in [z,G]$, and any $y_2, \dots, y_m \in \gamma_2(G)$, there exists
some $\alpha \in \Aut_c(G)$ such that
\begin{equation}
\label{equa}\alpha(x_i) = x_iy_i  \text{ for all } i = 1,2,\dots,m.
\end{equation}

Since the nilpotency class of $G$ is $\ge 3$, it follows from Proposition 
\ref{prop1a} that there exists a subgroup $H$ of index $p$ in $G$ such 
that $\gamma_{2}(H) = \gamma_{2}(G)$ and $\gamma_{2}(G) \not \leq \Z(H)$. 
The following three cases arise:

Case 1.  $z \in H$. Since $\gamma_{2}(G) = \gamma_{2}(H)$, $z \in
H-\gamma_{2}(G)$ can be extended to  a  minimal generating set $\{z = x_1,
x_2, \dots, x_m\}$ for $G$ such that $\gen{x_1,x_2, \dots, x_{m-1}} = H$ and
$x = x_m \notin H$. Since the group $G$ is almost Camina and $[z,G] \subset
\gamma_2(G)$, it follows that $[x,G] = \gamma_{2}(G)$. Pick any element 
$y \in \gamma_2(G) - \Z(H)$. Now we apply \eqref{equa}
with $y_1 = y_2 = \dots = y_{m-1} = 1$ and $y_m = y$ to get an automorphism
$\alpha \in \Aut_c(G)$ such that
\[\alpha(x_{i}) = x_{i},\;1 \leq i \leq m-1 \;\;\text{and} \;\; \alpha(x) =
xy.\]
Then $\alpha$ centralizes $H = \gen{x_1, \cdots, x_{m-1}}$, as well
as $G/H = \gen{xH}$.  But any automorphism of $G$ centralizing both
$G/H$ and $H$ must send the generator $x$ of $G$ modulo $H$ to $xw$, for
some element $w \in \Z(H)$ \cite[Satz I.17.1]{bH67}. Since $\alpha(x) = xy$ 
and $y \notin \Z(H)$, we get a contradiction.

Case 2.  $z \notin H$ and $[z,G] \not \subseteq \Z(H)$.  
Since $[z,G] \not\subseteq \Z(H)$, there exists a non-trivial element
$y \in [z,G] - \Z(H)$. We can also pick a minimal set of
generators $\{x_1, x_2,\dots, x_m\}$  for $G$ such that $x_1 = z$ and 
$\gen{x_2, x_3, \dots, x_m} = H$. Now we apply \eqref{equa} with $y_1 = y$ 
and $y_2 = y_3 = \dots = y_m = 1$ and  get an automorphism 
$\alpha \in \Aut_c(G)$ such that $\alpha(x_{i}) = x_{i},\;1 \leq i \leq m-1$
and $\alpha(x_1) = x_1y$.  Then again, as in Case 1, $\alpha$ centralizes 
$H = \gen{x_2, \cdots, x_{m}}$, as well as $G/H = \gen{x_1H}$ and therefore
$\alpha$  must send the generator $x_1$ of $G$ modulo $H$ to $x_1w$, for
some element $w \in \Z(H)$. Since $\alpha(x_1) =x_1y$ with $y \not\in \Z(H)$, 
we get a contraction in this case too.

Case 3. $z \notin H$ and $[z,G] \subseteq \Z(H)$. 
We can choose a minimal generating set $\{z=x_1, x_2, \dots, x_m\}$ for $G$
such that $H = \gen{x_2, \dots, x_m}$.
 We know that $[z,G] \subseteq \Z(H) \cap \gamma_2(G) < \gamma_2(G)$. 
Since $|[z,G]| = p^{m-1}$ and $|\gamma_2(G)| = p^m$, this forces 
$[z,G] = \Z(H) \cap \gamma_2(G) = \Z(H) \cap \gamma_2(H)$.  
More precisely, $[z,G] = \Z(H)$. For proving this, it is sufficient to prove
that $\Z(H) \leq \gamma_2(H)$. Suppose the contrary, i.e., $\Z(H) \not\leq
\gamma_2(H)$. Then there exists a non-trivial element $u \in
\Z(H)-\gamma_2(H)$. Since $[z,G] \subset \gamma_2(G)$ and $z \not\in H$, 
it follows that for this $u \in H-\gamma_2(H)$ we have $[u,G] =\gamma_2(G)$.
Therefore, since $z^p \in H$ and $u \in \Z(H)$,
\[p^3 \le p^m = |\gamma_2(G)| = |[u,G]| = |[u,H\gen{z}]| = |[u,\gen{z}]| \le p,\]
which is absurd. Hence $\Z(H) \leq \gamma_2(H)$.

If $C_{G}(z)\gamma_{2}(G) \neq \gen{z}\gamma_{2}(G)$, then there exists
an element $y \in H-\gamma_2(G)$ such that $yz=zy$. 
Since $[z,G] \subset \gamma_2(G)$ and $z \not\in H$, it follows from Lemma
\ref{acaminalemma} that $[y,G] = \gamma_2(G)$.
Since $y \in H$, $\Z(H) \le \C_H(y) \le \C_G(y)$. 
Also $y, z \in \C_G(y)$. So $\Z(H)\gen{y,z} \subseteq
\C_G(y)$. Thus $|\C_G(y)| \ge |\Z(H)\gen{y,z}| \ge p^{m+1}$. This implies that
\[|\gamma_2(G)| = |[y,G]| = |y^G| = |G|/|\C_G(y)| \le p^{m-1},\]
  which is a contradiction. So we can assume that 
$C_{G}(z)\gamma_{2}(G) = \gen{z}\gamma_{2}(G)$.
Since  $C_{G}(z)\gamma_{2}(G) = \gen{z}\gamma_{2}(G)$ and
$|\C_G(z)| = p^{m+1} = |\gen{z}\gamma_2(G)|$, it follows that
$\C_G(z) = \gen{z}\gamma_2(G)$. 
Now $G = \C_G(z)H$, so that $[z,H] = [z,G] = \Z(H)$. 
Set $[z,H] = \Z(H) = N$.  Since $z$
centralizes $\gamma_2(G)$, it follows that the map $\eta \colon
h\gamma_2(G) \mapsto [z,h]$ is a well defined epimorphism of $\bar H =
H/\gamma_2(G)$ onto $N = [z,H]$.  But both $\bar H$ and $N$ have order
$p^{m-1}$. Thus $\eta$ is an isomorphism, and $N$ is elementary abelian,
with the $\eta(x_i) = [z,x_i]$, for $i = 2,3,\dots, m$, as a basis.

The factor group $\gamma_2(G)/N = \gamma_2(H)/N$ has order $p$.
If $y \in H - \gamma_2(G)$, then $\gamma_2(G) = [y,G] = [y, \gen{z}H] \equiv
[y,H] \mod N$.  It follows that $H/N$ is an extra-special group of order
$p^m$.  So $m$ is odd and $m \ge 3$. Since $H/N = H/\Z(H)$ has class $2$,
it follows that $\gamma_3(H) \leq \Z(H)$. Hence the nilpotency class of $H$
is $3$ and hence from Lemma \ref{lemma2}, the class of $G$ is $3$.

Suppose that $m > 3$. Then we may suppose that $[x_2,x_3]$
generates $\gamma_2(H)$ modulo $N$, and that $[x_2, x_i] = [x_3,x_i] \in
N$ for all $i = 4,5,\dots, m$.  The identity
\[  [[x_2,x_3], x_i][[x_i, x_2], x_3][[x_3, x_i], x_2] = 1   \]
holds for any such $i$, since the class $H$ is $3$.  Here $[x_i,x_2]$ and
$[x_3,x_i]$ lie in $N = \Z(H)$, so that $[[x_i,x_2],x_3] = [[x_3,
x_i],x_2] = 1$.  Therefore $[[x_2,x_3], x_i] = 1$. Since $x_i$ certainly
centralizes $N$, this implies that $[\gamma_2(H), x_i] = 1$ for $i =
4,5,\dots,m$.  Since $m \ge 5$, we may assume (by re-indexing the generators, if
necessary) that $[x_4,x_5]$ also
generates $\gamma_2(H)$ modulo $N$.  Since $[x_i,x_j] \in N$ for $i =
2,3$ and $j = 4,5$, a similar argument shows that $[\gamma_2(H), x_2] =
[\gamma_2(H), x_3] = 1$. Thus $\gamma_2(H) \le \Z(\gen{x_2,x_3, \dots, x_m})
= \Z(H)$, which is impossible because $H$ has class $3$. Therefore $m = 3$.

Now $|N| = p^{m-1} = p^2$, and $|\gamma_2(H)| = p|N| = p^3$.  
It follows that $|H| = p^5$ and therefore $|G| = p^6$.  
The element $[x_2, x_3]$ must generate $\gamma_2(H)$ modulo $\gamma_3(H)$. 
Since the class of $H$ is $3$, this forces the elementary abelian group 
$\gamma_3(H) \le N$ to equal $N$, and to
have $[[x_2,x_3],x_2]$ and $[[x_2,x_3],x_3]$ as a basis.  Here we claim that
 $p$ is odd. 
Suppose that $p = 2$. If exponent of $H/\gamma_3(H)$ is $2$, then 
$H/\gamma_3(H)$ must be abelian. Therefore $\gamma_2(H) \le \gamma_3(H)$,
which is not possible, since the class of $H$ is $3$. 
Thus there is some element
$u \in H - \gamma_2(H)$ such that $u^2 \in \gamma_2(H) - \gamma_3(H)$.
So $\gamma_2(H)/\gamma_3(H) = \gen{u^2\gamma_3(H)}$. 
Since $\gamma_2(H) = \Phi(H)$, $u$ can be extended to a minimal generating set
$\{u,v\}$ for $H$. Now from what we have had just above, $[u,v]$ must generate
$\gamma_2(H)$ modulo $\gamma_3(H)$ and $[[u,v],u]$ must be non-trivial. 
Therefore $[u,v] \equiv u^2$ modulo $\gamma_3(H) (= N = \Z(H))$. 
Thus $[[u,v],u] = [u^2,u] = 1$, which is a contradiction. 
Hence our claim is true.
It now follows from Lemma \ref{lemma2} that $\gamma_3(G) = \gamma_3(H)$.
Since the class of $G$ is $3$ and $\Z(G) < \gamma_2(G)$, we must have 
$\gamma_3(G) = \Z(G)$. 

Thus if  $m > 2$, then $m = 3$ and $p>2$. Furthermore, $|G| = p^6$, and $G$
has a minimal set of generators $x_1, x_2, x_3$ satisfying
following properties:
\begin{subequations}
\begin{align}
& \text{$G = \gen{x_1, x_2, x_3}$.}\label{c1}\\
& \text{$H = \gen{x_2,x_3}$ is a maximal subgroup of $G$ such that 
$\gamma_2(G) = \gamma_2(H)$.}\label{c2}\\
&\text{$|\gamma_2(G)| = p^3$ and $|\gamma_3(G)| = p^2$.}\label{c3}\\
&\text{$[x_1,x_2] \neq 1$ and $[x_1,x_3] \neq 1$ generate 
$\gamma_3(G)$.}\label{c4}\\
& \text{$\beta' = [x_2,x_3]$ generates $\gamma_2(G)$ modulo 
$\gamma_3(G)$.}\label{c5}\\
& \text{$\beta_1^{\prime} = [\beta^{\prime}, x_2]$ and  
$\beta_2^{\prime} = [\beta^{\prime}, x_3]$ generates $\gamma_3(G)$;}
\label{c6}\\
& \text{$[z,G] = N = \Z(H) = \gamma_3(H) = \gamma_3(G) = \Z(G)$.}\label{c7}
\end{align}
\end{subequations}

Set $x_1 = \alpha$, $x_2 = \alpha_1$, $x_3 = \alpha_2$, $\beta' = \beta$,
$\beta_1' = \beta_1$ and $\beta_2' = \beta_2$.
Now it follows from \cite{rJ80} that there are only three isoclinism families
$\phi_{19}$, $\phi_{20}$ and $\phi_{21}$ such that any group $G_1$ from these
families satisfy the conditions \eqref{c1}, \eqref{c2} and \eqref{c3}.
But $\phi_{21}$ is the only isoclinism family such that the groups $G_1$
from it satisfy  conditions \eqref{c1}--\eqref{c4}. 
Thus our groups $G$, if they exist, must lie in $\phi_{21}$.
Hence $G$ is isoclinic to $R$.      \hfill $\Box$

\end{proof}

A finite group $G$ is said to be \emph{purely non-abelian} if it does not 
have a non-trivial abelian direct factor. An automorphism $\phi$ of a group 
$G$ is called
\emph{central} if $g^{-1}\phi(g) \in \Z(G)$ for all $g \in G$. The set of all
central automorphisms of $G$, denoted by $\Autcent(G)$, is a normal subgroup
of $\Aut(G)$.

\begin{prop}\label{prop2b}
Let $G$ be any group of order $p^6$ which is isoclinic to $R$.
Then $G$ is an almost Camina group and equality holds in \eqref{gbineq} 
for $G$.
\end{prop}
\begin{proof}
Since $|G| = p^6$ and $G$ is isoclinic to $R$, $G$ belongs to the isoclinism
family $\phi_{21}$ of \cite{rJ80}. 
Now it is easy to check that $G$ satisfies all the
conditions \eqref{c1}--\eqref{c7} with the setting $x_1 = \alpha$, 
$x_2 = \alpha_1$, $x_3 = \alpha_2$, $\beta' = \beta$,
$\beta_1' = \beta_1$ and $\beta_2' = \beta_2$.    
Also notice that the groups $\gamma_{3}(G)$ and $G/\gamma_2(G)$ are 
elementary abelian of order $p^2$ and $p^3$ respectively. 
It follows from \cite{rJ80} that there is no $g \in G$ such that $|g^G| = p$.
Thus for any $g \in G-\Z(G)$, $|g^G| \ge p^2$. 

Let $g \in \gen{\alpha}\gamma_2(G)-\gamma_2(G)$.
Then $g = \alpha^iu$, where $u \in \gamma_2(G)$ and $1 \leq i \leq p-1$. Let
$g'$ be an arbitrary element of $[g,G] = [\alpha^iu,G]$. 
Then $g' = [\alpha^iu,v]$ for some $v \in G$. 
Therefore
\[g' = [\alpha^iu,v] = [\alpha^i,v]^u[u,v] = [\alpha,v]^i[u,v] \in \Z(G),\]
since $[\alpha,v] \in [\alpha,G]$, $[u,v] \in \gamma_3(G)$ and 
from \eqref{c7} we know that $[\alpha,G] = \Z(G) = \gamma_3(G)$ 
(with above setting).  Thus $[g,G] \subseteq \Z(G)$. 
Since $g \in G-\gamma_2(G) \subset G-\Z(G)$, $|[g,G]| \ge p^2$. 
Hence $[g,G] = \Z(G)$, since $|\Z(G)| = p^2$. Thus for any
$g \in \gen{\alpha}\gamma_2(G) - \Z(G)$, we have that $|g^G| = p^2$.
Since $|\gen{\alpha}\gamma_2(G) - \Z(G)| = p^4-p^2 = p^2(p^2-1)$, it follows
that $G$ has at least $p^2-1$ conjugacy classes of length $p^2$.
But from \cite[$\phi_{21}$]{rJ80} it follows that $G$ has $p^2$ conjugacy
classes  of length $1$, $p^2-1$ classes of length $p^2$ and $p^3-p$ 
classes of length $p^3$. 
Since $|\Z(G)| = p^2$, therefore $|g^G| = p^3$ for all $g \in G-
\gen{\alpha}\gamma_2(G)$. Hence it follows that $G$ is an almost Camina group.

Since $\Z(G) \subseteq \gamma_2(G)$, $G$ is purely non-abelian. Then it
follows from \cite[Theorem 1]{AY65} that 
$|\Autcent(G)| = |\Hom(G/\gamma_2(G), \Z(G))|$.
Since both $G/\gamma_2(G)$ and $\Z(G)$ are elementary, we have from
Lemma \ref{hlemma} that $|\Autcent(G)| = p^6$. 
We claim that $\Autcent(G) \subseteq \Aut_c(G)$.
Let $\phi \in \Autcent(G)$. Then $\phi(g) = g$ for all $g \in \gamma_2(G)$.
So let $g \in G-\gamma_2(G)$. Since $[g,G] = \Z(G)$ for all $g \in
\gen{\alpha}\gamma_2(G) - \gamma_2(G)$ and $[g,G] = \gamma_2(G)$ for all $g \in
G-\gen{\alpha}\gamma_2(G)$, we have that $\Z(G) \le [g,G]$ for all 
$g \in G-\gamma_2(G)$.  
Since $\phi \in \Autcent(G)$,  $g^{-1}\phi(g) \in \Z(G) \le [g,G]$.
Thus $\phi(g) \in g^G$ for all $g \in G$. This proves that $\phi \in
\Aut_c(G)$ and our claim follows. 

Let $i_2$ and $i_3$ be the inner automorphisms of $G$ induced by $x_2$ and
$x_3$ respectively. 
Since $x_2^{-1}i_3(x_2) = [x_2,x_3]$ and 
$x_3^{-1}i_2(x_3) = [x_3,x_2]$ lie outside $\Z(G)$, it follows that $i_j$
is not central for $j = 2, 3$.
If $i_2^k = i_3$, $1 \le k < p$, then $i_2^k(x_2) = i_3(x_2)$. 
This implies that $x_2 = x_3^{-1}x_2x_3$. So $[x_2,x_3]=1$, which is not true.
Therefore $i_2^k$ and $i_3$ are distinct, $1 \le k < p$. 
Similarly $i_3^k$ and $i_2$ are distinct, $1 \le k < p$.
Let $\I = \gen{i_2, i_3}$. Then $\I \le \Inn(G)$ such that $|\I| \ge p^2$.  
Now let $i = i_2^{k_2}i_3^{k_3} \in \I$, where $0 \le k_2 < p$ and 
$1 \le k_3 < p$. If $i \in \Autcent(G)$, then $x_2^{-1}i(x_2) \in \Z(G)$.
So $[x_2, x_3^{k_3}x_2^{k_2}] \in \Z(G)$. 
This implies that $[x_2, x_3^{k_3}]^{x_2^{k_2}} \in \Z(G)$. 
Which gives that $[x_2, x_3^{k_3}] \in \Z(G)$ ($= \gamma_3(G)$) and 
therefore $[x_2, x_3] \in \Z(G)$, since
$1 \le k_3 < p$. But this is not true. Thus $i$ is not central. This proves
that $|\I \Autcent(G)/\Autcent(G)| \ge p^2$.
Thus  $|\Aut_c(G)| \ge |\Autcent(G) \I| \ge p^8$. But from \eqref{bineq},
it follows that $|\Aut_c(G)| \le p^8$. 
Hence $|\Aut_c(G)| = p^8$ and equality hold in 
\eqref{gbineq}. \hfill $\Box$

\end{proof}

\begin{lemma}\label{clemma}
The group $W$ is a Camina special $p$-group of order $p^6$.
\end{lemma}
\begin{proof} Notice that $W$ is a special $p$-group such that $|W| = p^6$
and $|\Z(W)| = p^2$. From \cite{wB13} we have
\[|\Aut_{c}(W)| = p^8 = p^{m(n-m)},\]
for $m = 2$ and $n = 6$. Thus equality holds in \eqref{lbineq} and $W$ is
a Camina group by Theorem \ref{thm2}. \hfill $\Box$

\end{proof}

Let $\Cl(G)$ be the set of all conjugacy classes of $G$. Let
\[\cl(G) = \{\cl(x)\;|\;x \in G\} = \{1,\;p^{n_{1}},\; \cdots, \;p^{n_{r}}\},\]
where $\cl(x) = |x^G|$ and $0 < n_{1} < \cdots < n_{r}$. In this situation 
$G$ is called a group of \emph{(conjugate) type} 
$\{1,\;p^{n_{1}},\; \cdots, \;p^{n_{r}}\}$.

\begin{lemma}\label{cslemma}
Let $G$ be any special $p$-group isoclinic to $W$. Then $G$ is a Camina
special $p$-group such that $|G|= p^6$ and $|\gamma_{2}(G)| = p^2$.
\end{lemma}
\begin{proof} Since $G$ is special, $\gamma_{2}(G) = \Z(G) = \Phi(G)$. 
And since $G$ is isoclinic to $W$, $\gamma_{2}(G)\cong \gamma_{2}(W)$ and
$G/\Z(G) \cong W/\Z(W)$. 
It follows from \cite{pH40} that $G$ and $W$ have the same conjugate type. 
Thus $G$ is of conjugate type $\{1, |\gamma_{2}(G)|\}$.
This implies that $G$ is a Camina group.
Now $|\gamma_{2}(G)| = |\gamma_{2}(W)| = p^2$ and $|G/\gamma_{2}(G)| 
= |W/\gamma_{2}(W)| = p^4$, since both $G$ and $W$ are special $p$-groups.
Hence
\[|G| = |G/\gamma_{2}(G)||\gamma_{2}(G)| = p^6.\] \hfill $\Box$

\end{proof}

\begin{lemma} \label{goodlemma}
Up to isoclinism there is only one Camina special $p$-group $H$ of order 
$p^6$ such that $|Z(H)| =p^2$. 
\end{lemma}
\begin{proof}

Let $H$ be any Camina special $p$-group such that 
$|H|=p^6$ and $|Z(H)| =p^2$. 
Since isoclinism is an equivalence relation and the group $W$  
is a Camina special $p$-group (Lemma \ref{clemma}), 
it is sufficient to prove that $H$ is isoclinic to $W$.
Since $H$ is a Camina special group, $|Z(H)| = |\gamma_{2}(H)|
= |\phi(H)| = p^2$. Also $|Z(W)| = |\gamma_{2}(W)| = |\phi(W)| = p^2$. 
$|W/Z(W)| = |H/Z(H)| = p^4$ and  both of these factor groups are elementary. 
Since there is only one elementary group of given order upto isomorphism, it 
follows that $\gamma_{2}(W) \cong \gamma_{2}(H)$ and $W/Z(W) \cong H/Z(H)$.

Set $\bar H = H/\Z(H)$, and consider the bilinear map 
$a = a_H \colon \bar H \times \bar H \to [H,H] =
\Z(H)$ induced by commutation.  The Camina condition is that $a(\bar h,
\bar H) = \Z(H)$ for any $\bar h \in \bar H^{\#} = \bar H -
\{1\}$. Because $\bar H$ and $\Z(H)$ are elementary abelian $p$-groups
with ranks $4$ and $2$, respectively, it follows that 
\[ \bar h^{\perp} = \{\, \bar h' \in \bar H \mid a(\bar h, \bar h') =
1\,\} \]
is elementary abelian of rank $2$, for any $\bar h \in \bar H^{\#}$.
We claim that $a(\bar h^{\perp}, \bar h^{\perp}) = 1$, so that 
$\bar h^{\perp} = (\bar h')^{\perp}$, 
for any non-trivial $\bar h' \in \bar h^{\perp}$.
For, $\bar h^{\perp}$ has rank
$2$, it has a basis of the form $\bar h, \bar h''$, where $\bar h'' \in
\bar h^{\perp} - \gen{\bar h}$.  Then $a(\bar h, \bar h'') = 1 = a(\bar h'',
\bar h)$, since $\bar h'' \in \bar h^{\perp}$.  Also $a(\bar h, \bar h)
= 1 = a(\bar h'',\bar h'')$. Since $\bar h$ and $\bar h''$ generate $\bar
h^{\perp}$, this and the bilinearity of $a$ imply that $a(\bar
h^{\perp},\bar h^{\perp}) = 1$.  It follows that $\bar h^{\perp}
\subseteq (\bar h')^{\perp}$ for any non-trivial $\bar h' \in \bar
h^{\perp}$. Since both $\bar h^{\perp}$ and $(\bar h')^{\perp}$ have
order $p^2$, this implies the equality $(\bar h')^{\perp} = \bar
h^{\perp}$.
Thus we have an equivalence relation on $\bar H^{\#}$, whereby two
elements $\bar h, \bar h' \in \bar H^{\#}$ are equivalent if and only if
$a(\bar h, \bar h') = 1$.  The equivalence class of $\bar h \in \bar
H^{\#}$ is just the set $(\bar h^{\perp})^{\#}$ of non-trivial elements
in $\bar h^{\perp}$.

Take any $y_1 \in H - \Z(H)$, with image $\bar y_1 \in \bar H^{\#}$.  
Then take any $y_2 \in H$ whose image $\bar y_2 \in \bar H$ completes 
a basis $\bar y_1, \bar y_2$ for $\bar y_1^{\perp}$.  
For $y_3$ we take any element in $H$
whose image $\bar y_3$ lies in $\bar H - \bar y_1^{\perp}$. For $y_4$ we
take any element of $H$ whose image $\bar y_4$ completes a basis $\bar
y_3, \bar y_4$ for $\bar y_3^{\perp}$. The intersection $\bar
y_1^{\perp} \cap \bar y_3^{\perp}$ is trivial, since $\bar y_3 \notin
\bar y_1^{\perp}$. By counting ranks, we see that $\bar H$ must be the
direct product $\bar y_1^{\perp} \times \bar y_3^{\perp} = \gen{\bar y_1,
\bar y_2} \times \gen{\bar y_3, \bar y_4}$. 
Thus it follows that  $H = \gen{y_{1}, \dots, y_{4}}$ such that 
$[y_{1},\;y_{2}] = 1$ and $[y_{3},\;y_{4}] = 1$. 
Now
\[\gamma_{2}(H) = \gen{[y_{i},\;y_{j}]^H | 1 \leq i < j \leq 4} = 
\gen{[y_{i},\;y_{j}] | 1 \leq i < j \leq 4},\]
since $\gamma_{2}(H) = Z(H)$. We know that 
$|C_{H}(y_{i})| = |H|/|y_{i}^H| = p^4$. 
Thus for any $y_{i}$ there is exactly one $y_{j}$, $i \neq j$ such that
$[y_{i}, y_j] = 1$. Which implies that all the elements $[y_{1}, y_3]$,
$[y_{1}, y_4]$, $[y_{2}, y_3]$ and $[y_{2}, y_4]$ are non-trivial.

Since $\bar W$ and $\bar H$ are elementary, we can define an
isomorphism  $\phi : \bar W  \to  \bar H$ such that 
$\phi (\bar x_{i}) = \bar y_{i}$, $1 \le i \le 4$.
Now define a map $\theta : [W, W] \to [H,H]$
such that $\theta ([x_{i},x_{j}]) = [y_{i}, y_{j}]$, $1 \leq i < j \leq 4$.
To prove that $\theta$ is an isomorphism, it is sufficient to show that
$\theta$ is well defined.
To prove that $\theta$ is well defined we must show that the elements 
$[x_i,x_j]$
and the elements $[y_i,y_j]$ satisfy the same relations for all $i, j$ such
that $1 \leq i < j \leq 4$. The rest of the proof is devoted to this.

Write both $\bar H = H/[H,H]$ and
$[H,H]$ additively, obtaining vector spaces $U = \bar H^+$ and $V =
[H,H]^+$ of dimensions $4$ and $2$, respectively, over the field
$\mathbb Z_p$ of $p$ elements.  Then commutation in $H$
induces a strongly alternating, $\mathbb Z_p$-bilinear map $a$ of $U
\times U$ into $V$.  The vector space $U$ is the direct sum $U' \dotplus
U''$ of two-dimensional subspaces $U' = \gen{\bar y_1, \bar y_2 }^+$ and
$U'' = \gen{\bar y_3, \bar y_4 }^+$, which satisfy
\begin{equation}\label{e1}
 a(U',U') = a(U'',U'') = 0.
\end{equation}
Furthermore, $a(u',u'') \neq 0$ for all non-zero $u'
\in U'$ and $u'' \in U''$. Since both $U''$ and $V$ have dimension $2$,
this implies that the map 
\[ T_{u'} \colon u'' \mapsto T_{u'}(u'') = a(u',u'') \] 
is an isomorphism of the vector space $U''$ onto $V$, for
any $u' \neq 0$ in $U'$. It follows that any ordered basis $u'_1,u'_2$
for $U'$ determines an automorphism 
\[ S = S_{u'_1,u'_2} = T_{u'_1}^{-1}T_{u'_2} \]
of the vector space $U''$.

        We claim that $S$ has no eigenvalue $n \in \mathbb Z_p$.
Suppose that such an $n$ exists.  Then there is some non-zero
eigenvector $u'' \in U''$ such that $S(u'') = nu''$.  This implies that
\[ a(u'_2, u'') = T_{u'_2}(u'') = T_{u'_1}S(u'') = T_{u'_1}(nu'') =
a(u'_1,nu'') = a(nu'_1,u'') \in V. \] 
So $-nu'_1 + u'_2 \in U'$ and $u'' \in U''$ are non-zero elements such
that $a(-nu'_1 + u'_2, u'') = 0$. This is impossible by
hypothesis. Therefore our claim is correct.

Now $S$ is a linear transformation of a vector space $U''$ of
dimension $2$ over $\mathbb Z_p$. Furthermore, $S$ has no eigenvalue in
$\mathbb Z_p$. It follows that the $\mathbb Z_p$-algebra $\mathbb
Z_p[S]$ of linear transformations of $U''$ generated by $S$ is
isomorphic to the field $\mathbb F_{p^2}$ of $p^2$ elements.  Fix an
isomorphism $i$ of $\mathbb F_{p^2}$ onto $\mathbb Z_p[S]$, and a
generator $f$ for $\mathbb F_{p^2}$ over $\mathbb Z_p$.  Then $1$ and
$f$ form a $\mathbb Z_p$-basis for $\mathbb F_{p^2}$, so that $i(1) = I$
and $i(f)$ form a $\mathbb Z_p$-basis for $\mathbb Z_p[S]$. Since $I$
and $S$ also form a $\mathbb Z_p$-basis for $\mathbb Z_p[S]$, there are
some elements $k,m \in \mathbb Z_p$ such that $i(f) = kI + mS$ and $m
\neq 0$.

Evidently $u'_1$ and $u'_3 = ku'_1 + mu'_2$ also form an ordered
basis for $U'$ over $\mathbb Z_p$. The $\mathbb Z_p$-bilinearity of $a$
implies that 
\[ T_{u'_3}(u'') = a(u'_3, u'') = ka(u'_1,u'') + ma(u'_2,u'') =
(kT_{u'_1} + mT_{u'_2})(u'') \in V \]
for any $u'' \in U''$. Therefore $T_{u'_3} = kT_{u'_1} + mT_{u'_2}$, and
\[ S_{u'_1,u'_3} = T_{u'_1}^{-1}T_{u'_3} = T_{u'_1}^{-1}(kT_{u'_1} +
mT_{u'_2}) = kI + mS = i(f). \] 
Thus we can replace our original ordered basis $u'_1,u'_2$ for $U'$ by
$u'_1,u'_3$, and assume from now on that
\[ S = i(f).  \]

Let $X^2 + bX + c \in \mathbb Z_p[X]$ be the minimal polynomial
of $f$ over $\mathbb Z_p$. Then $X^2 + bX + c$ is both the minimal and
the characteristic polynomial of $S$ over $\mathbb Z_p$.  In particular,
$S^2 = -cI - bS$.  Let $u''_1$ be any non-zero element of $U''$, and
$u''_2$ be $S(u'_1)$. Then $u''_1$ and $u''_2$ form a $\mathbb
Z_p$-basis for $U''$.  Set 
\begin{equation}\label{e2}
 v_1 = a(u'_1,u''_1) = T_{u'_1}(u''_1) \quad \text{and} \quad v_2 =
a(u'_1, u''_2) = T_{u'_1}(u''_2).
\end{equation}
Then $v_1$ and $v_2$ form a ${\mathbb Z_p}$-basis for $V$. Furthermore, 
\begin{equation}\label{e3}
 a(u'_2, u''_1) = T_{u'_2}(u''_1) =
T_{u'_1}T_{u'_1}^{-1}T_{u'_2}(u''_1)
=  T_{u'_1}S(u''_1) =
T_{u'_1}(u''_2) = v_2
\end{equation}
and 
\begin{eqnarray}
 a(u'_2,u''_2) &=& T_{u'_2}(u''_2) =
T_{u'_1}T_{u'_1}^{-1}T_{u'_2}S(u''_1) = T_{u'_1}S^2(u''_1)\label{e4}\\
&=& T_{u'_1}(-cI -bS)(u''_1) = -cv_1 -bv_2.\nonumber
\end{eqnarray}
Now we have constructed a $\mathbb Z_p$-basis $u'_1, u'_2,
u''_1,u''_2$ for the vector space $U = \bar H^+$, and a $\mathbb
Z_p$-basis $v_1,v_2$ for $V = [H,H]^+$, such that the values of the
strongly alternating $\mathbb Z_p$-bilinear function $a \colon U \times
U \to V$ are completely determined by the equations \eqref{e1} -- \eqref{e4}. 
Since these equations do not depend on the group $H$, it follows that
the elements $[x_i,x_j]$ and the elements $[y_i,y_j]$ satisfy the same 
relations for all $i, j$ such that $1 \leq i < j \leq 4$. Hence $\theta$ is 
well defined and therefore an isomorphism of $[W,W]$ onto $[H,H]$. 
Thus it follows that $(\phi, \theta)$ is an isoclinism of $W$ onto $H$. 
This completes the proof of the lemma.   
\hfill $\Box$

\end{proof}

Now we prove our main theorem `Theorem B'.
\begin{thm}\label{thm6}
Let $G$ be a non-abelian finite $p$-group of order $p^n$. Then
equality holds in \eqref{gbineq} if and only if one of the following holds: 
\begin{subequations}
\begin{align}
& \text{$G$ is an extra-special $p$-group of order $p^3$;}\\
& \text{$G$ is  a group of nilpotency class $3$ and order $p^4$;}\\
& \text{$G$ is a Camina special $p$-group isoclinic to the group $W$ and
    $|G|=p^6$;}\\
& \text{$G$ is isoclinic to $R$ and $|G|=p^6$.}
\end{align}
\end{subequations}
\end{thm}
\begin{proof} Let $G$ be an extra-special group of order $p^3$.
Then 
\[|\Aut_{c}(G)| = |\Inn(G)| = p^2 = p^{\frac{3^2 - 1}{4}}.\]
 Hence equality holds in \eqref{gbineq}. 
Now suppose that $G$
is a group of order $p^4$ and of class $3$. Then $|\Z(G)| = p$. Notice that
$\Aut_{c}(G) = \Inn(G)$. Thus
\[|\Aut_{c}(G)| = |\Inn(G)| = |G|/|\Z(G)| = p^3 = p^{\frac{n^2-4}{4}},\]
where $n = 4$. Hence equality holds in \eqref{gbineq}.
Next assume that $G$ is a Camina special $p$-group isoclinic to the group 
$W$.
Then it follows from Lemma \ref{clemma} and Lemma \ref{cslemma} that  
$|G| = p^6$ and $|\gamma_{2}(G)| = p^2$.
By Theorem \ref{thm3} we have
\[|\Aut_{c}(G)| = p^{2(6-2)} = p^8 = p^{\frac{n^{2} - 4}{4}},\]
where $n = 6$. So equality holds in \eqref{gbineq}. 
Finally suppose that $G$ is isoclinic to $R$. Then it follows
from Proposition \ref{prop2b} that equality holds in \eqref{gbineq}
for this group $G$.

Conversely suppose that equality holds in \eqref{gbineq}. When $n$ is odd,
equality in \eqref{gbineq} implies that equality holds in \eqref{lbineq}.
But when $n$ is even equality in \eqref{gbineq} may or may not imply equality
in \eqref{lbineq}. We consider two separate cases: 
(i) equality holds in \eqref{lbineq};
(ii) equality does not hold in \eqref{lbineq}.

First we consider the case (i). Since equality holds in \eqref{lbineq}, 
it follows from 
Theorem \ref{thm2} that $G$ is a Camina special $p$-group. 
Using  Theorem 3.2 of \cite{iM81} we have that $n-m$  is even and 
$n-m \geq 2m$, where $|G| = p^n$ and $|\gamma_{2}(G)| = p^m$. 
Thus we get inequality \eqref{goodineq}. 

Equality  in \eqref{gbineq} gives us
\begin{equation}\label{gbeq}
|\Aut_{c}(G)| =
   \begin{cases}
     p^{\frac{n^{2}-4}{4}},  &\text{if $n$ is even;}\\
     p^{\frac{n^{2}-1}{4}}, &\text{if $n$ is odd.}
   \end{cases}
\end{equation}
First suppose that $n$ is even. Then \eqref{gbeq} holds only when 
$m = (n-2)/2$ or $m = (n+2)/2$. The value $m = (n+2)/2$ clearly contradicts 
\eqref{goodineq}. Thus only possibility is $m = (n-2)/2$. Which, along with
\eqref{goodineq}, gives $n \leq 6$. Since $G$ is non-abelian, the only
possibilities are $n = 4,\; 6$. When $n = 4$, there exists no $m$ such that
$n-m$ is even and \eqref{goodineq} is satisfied. Therefore we are left only 
with $n = 6$.
Assume that $|G| = p^6$. Then  $|Z(G)| = |\gamma_{2}(G)| = p^2$.
Now it follows from Lemma \ref{goodlemma} that $G$ is isoclinic to  
the group $W$.

Next suppose that $n$ is odd. Then \eqref{gbeq} holds  only when 
$m = (n-1)/2$ or $m = (n+1)/2$. The choice
$m = (n+1)/2$ again contradicts \eqref{goodineq}. So only possibility is
$m = (n-1)/2$. Which, using \eqref{goodineq}, gives $n \leq 3$. Since $G$ is 
non-abelian, $n$ must be $3$. Since every non-abelian group of order $p^3$ is
extra-special, we are done in this case. 

Now we consider the case (ii). Since equality in \eqref{gbineq} implies
equality in \eqref{lbineq} when $n$ is odd, we have that $n$ must be even and 
$|\gamma_{2}(G)| = p^{\frac{n}{2}}$. So $|G|= p^{2m}$ and 
$|\gamma_{2}(G)| = p^{m}$, where $m = n/2 \ge 2$.
Now it follows from Lemma \ref{ecd1} that $G$ is an almost Camina group.
It then follows from Proposition \ref{prop2a} that either $G$ is a group of 
nilpotency class $3$ and order $p^4$ or $G$ is isoclinic to $R$. 
This completes the proof of the theorem.             \hfill $\Box$     

\end{proof}

\end{document}